\documentclass[11pt]{amsart}
\usepackage{amsmath, amsthm, amssymb, latexsym,url}
\usepackage{hyperref}

\allowdisplaybreaks

\author{Joseph Vandehey}
\thanks{Email: \href{mailto:vandehey.1@osu.edu}{\nolinkurl{vandehey.1@osu.edu}}\\ The Ohio State University}
\title[Uncanny normal numbers]{Uncanny subsequence selections that generate normal numbers}
\date{\today}

\newtheorem{thm}{Theorem}[section]

\theoremstyle{remark}

\begin{document}

\begin{abstract}
Given a real number $0.a_1a_2 a_3\dots$ that is normal to base $b$, we examine increasing sequences $n_i$ so that the number $0.a_{n_1}a_{n_2}a_{n_3}\dots$ are normal to base $b$. Classically it is known that if the $n_i$ form an arithmetic progression then this will work. We give several more constructions, including $n_i$ that are recursively defined based on the digits $a_i$.

Of particular interest, we show that if a number is normal to base $b$, then removing all the digits from its expansion which equal $(b-1)$ leaves a base-$(b-1)$ expansion that is normal to base $(b-1)$.
\end{abstract}

\maketitle

\section{Introduction}

A number $x\in  [0,1)$ with base $10$ expansion $x=0.a_1a_2a_3\cdots$ is said to be normal (to base $10$) if for any finite string $s=[c_1,c_2,\dots,c_k]$ we have that 
\[
\lim_{n\to \infty} \frac{\#\{0\le i \le n-1: a_{i+j} = c_j, 1\le j \le k\}}{n} = \frac{1}{10^k}.
\] More colloquially, a number is normal to base $10$ if every string of base $10$ digits appears with the expected frequency if the digits were chosen at random. 
Although almost all real numbers are normal to base $10$, we still do not know of a single commonly used mathematical constant, such as $\pi$, $e$, or $\sqrt{2}$, that is normal. A similar definition of normality holds for other base-$b$ expansions besides decimal.

By a classical result of Wall \cite{Wall}, we have that if $0.a_1a_2a_3\cdots$ is normal, then, for any integers $k,m\ge 1$,  the number $0.a_ka_{m+k}a_{2m+k}a_{3m+k}\cdots$ is also normal. In other words, normality is preserved along arithmetic progressions. The topic of interest in this paper will be that of \emph{normality-preserving selection rules}. Given a number $0.a_1a_2a_3\cdots$, a selection rule is a method of generating an increasing sequence $\{n_i\}_{i=1}^\infty$ (which may be dependent on the original number) and a corresponding real number $0.a_{n_1}a_{n_2}a_{n_3}\cdots$.  A selection rule is said to be normality preserving if $0.a_{n_1}a_{n_2}a_{n_3}\cdots$ is  normal whenever $0.a_1a_2a_3\cdots$ is normal.

The most common normality-preserving selection rules that have been studied are prefix selection rules. In prefix selection rules, a set $\mathcal{A}$ of finite-length base-$b$ strings is given, and then a positive integer $n$ is included in the sequence of $n_i$'s if and only if the string $[a_1,a_2,\dots,a_{n-1}]$ of the first $n-1$ digits belongs to the set $\mathcal{A}$. The arithmetic progression construction mentioned above is an example of such a prefix selection rule: it consists of \emph{all} strings of lengths $k-1$, $m+k-1$, $2m+k-1$, and so on. 

Agafonov \cite{Agafonov} showed that sets $\mathcal{A}$ that were generated in a particular way by finite automata would produce normality-preserving prefix selection rules. In fact Agafonov showed that a number is normal if and only if every number generated by such a selection rule would also be normal. This could be rephrased to say any $\mathcal{A}$ which is a \emph{regular language} would give a normality-preserving prefix selection rule. Merkle and Reimann \cite{MR} showed that there exist deterministic one-counter languages $\mathcal{A}$ and linear languages $\mathcal{A}$ which do not preserve normality. Kamae  \cite{Kamae} gave a complete characterizations of all prefix seelction rules which consist of all strings with a given set of lengths (so the sequence $n_i$ is the same for all normal numbers). Kamae and Weiss \cite{KW} gave a very general class of prefix selection rules which includes Agafonov's result as a special case.

Becher and Carton \cite{becher2013normality} study such subsequences of normal numbers from the perspective of automatas and compressibility. They give several results, including another generalization of Agafonov's work on prefix selection, in addition to showing that certain types of suffix selection also preserve normality. (Here, suffix selection means that $n$ is included in the sequence of $n_i$'s if and only if the sequence $[a_{n+1},a_{n+2},a_{n+3},\dots]$ belongs to a certain set $\mathcal{B}$.) On the other hand, they show that if $n$ is included in the sequence of $n_i$'s if and only if $a_{n-1}$ and $a_{n+1}$ both equal $0$, then this is not normality preserving.

One thing all of the above results have in common is that the choice of whether to include a given integer $n$ in the sequence of $n_i$'s does not depend on the value of $a_n$ itself. Such selection rules seem to be beyond the techniques of Kamae and Weiss, and if they are accessible to the techniques of Becher and Carton, this is not mentioned in the paper. 

In this paper, we will make use of some of the author's recent work on augmented dynamical systems to show that some such uncanny selection rules do preserve normality. These augmented systems are a skew-product over the original dynamical system, that build in a finite automata and allow the dynamical system to maintain a small amount of memory. This is similar in method to how Jager and Liardet proved Moeckel's Theorem \cite{JL}. We will use this technique to provide new and simple proofs of some results obtained by the above authors, and also some proofs of new normality-preserving selection rules. 

The first result we will show is not a new one (it can be seen to be a special case of the result of Kamae and Weiss); however, we include it in part because the proof is very simple and will be a template for future proofs.

\begin{thm}\label{thm:leap}
Let $x=0.a_1a_2a_3\cdots$ be normal to base $b$. Let $n_i$ be an increasing sequence recursively defined starting from an arbitrary $n_1 \in \mathbb{N}$ so that 
\[
n_{i+1} = n_i+1+a_{n_i}, \quad i\ge 1.
\]
Then the number $0.a_{n_1}a_{n_2}a_{n_3}\cdots$ is normal to base $b$.
\end{thm}

The next result appears to be new in the literature. It could be proven using a careful counting argument. We include it here as another demonstration of the method of proof. In fact, we will simply cite the proof of Theorem \ref{thm:leap} and only mention some minor variations.

\begin{thm}\label{thm:reduce}
Let $x=0.a_1a_2a_3\cdots$ be normal to base $b$. Let $n_i$ be the increasing sequence of positive integers defined so that $n=n_i$ for some $i$ if and only if $0\le a_n <b-1$. Then the number $0.a_{n_1}a_{n_2}a_{n_3}\cdots$ is normal to base $(b-1)$.
\end{thm}

One case of the above theorem is that if one starts with a base-$10$ normal number and removes all the $9$'s, one is left with a base-$9$ normal number. We leave the following as an open question: can every base-$9$ normal number be obtained in this fashion?

The last result of this paper, given below, appears not only to be new, but unprovable using any other methods so far known.

\begin{thm}\label{thm:modulo}
Let $x=0.a_1a_2a_3\dots$ be normal to base $b$. Let $L$ and $N$ be integers with $N\ge 2$. Let $n_i$ be the increasing sequence of positive integers defined so that $n=n_i$ for some $i$ if and only if 
\[
\sum_{j=1}^{n} a_j \equiv L \pmod{N} \qquad i\ge 1.
\]
Then the number $0.a_{n_1}a_{n_2}a_{n_3}\dots$ is normal to base $b$.
\end{thm}

\section{An augmented system}

We will require a result from a previous paper of the author \cite{ratmultCF}.

Let $T=T_b$ be the usual base-$b$ forward iteration map acting on the set $\Omega =[0,1)$. So $Tx \equiv b x\pmod{1}$.

Given a base $b$, we define the cylinder set $C_s$ for a string $s=[c_1,c_2,\dots,c_k]$ of base $b$ digits to be the set of points $x\in [0,1)$ with an expansion of the form $x=0.c_1c_2c_3\dots c_k a_{k+1}a_{k+2}\dots$. If $s$ is the empty string, then $C_s=[0,1)$. Cylinder sets will always be intevals, and we generally assume that the left endpoint is included in the interval, but not the right endpoint. It is easy to see that $\lambda(C_s)=b^{-k}$.

With these definitions (and a number $x=0.a_1a_2a_3\cdots$), we see that the term $\#\{0\le i \le n: a_{i+j} = c_j, 1\le j \le k\}$ in the definition of normality is equivalent to $\#\{0 \le i \le n : T^i x \in C_s\}$.  

The base-$b$ map $T$ on the set $\Omega$ is a typical example of an ergodic system---that is, the only measurable sets $E\subset \Omega$ such that $T^{-1} E = E$ have either full measure or null measure. This map also preserves the Lebesgue measure $\lambda$---that is, $\lambda(T^{-1}E) = \lambda(E)$ for all measurable $E\subset\Omega$. A consequence of the pointwise ergodic theorem is that for almost all $x\in \Omega$, we have that
\[
\lim_{n\to \infty} \frac{\#\{0 \le i \le n-1 : T^i x \in C_s\}}{n} = \lambda(C_s).
\]
In other words, almost all points $x$ are normal to base $b$.

We wish to extend the map $T$ to a skew-product transformation $\widetilde{T}$ on a larger domain $\widetilde{\Omega}=\Omega \times \mathcal{M}$, where $\mathcal{M}$ is some finite set. For any $(x,M)\in \widetilde{\Omega}$, we define \[\widetilde{T}(x,M) = (Tx,f(a_1(x),M)),\] where $a_1(x)$ is the first base-$b$ digit of $x$ and $f$ is a function from $\{0,1,2,\dots,b-1\}\times \mathcal{M}$ to $ \mathcal{M}$.  We also have a probability measure $\mu$ on $\widetilde{\Omega}$ that will be defined as being the product of the lebesgue measure on $\Omega$ and some weighted counting measure on $\mathcal{M}$. For easier readability, we will use $(E,M)$ to denote $E \times \{M\}$ for any measurable set $E\subset \Omega$, with measurability being determined by Lebesgue measure. 

We adapt our definition of normality on this space. We will say that $(x,M)\in \widetilde{\Omega}$ is $\widetilde{T}$-normal with respect to a measure $\mu$ on $\widetilde{\Omega}$, if for any cylinder set $(C_s,M')$ we have
\[
\lim_{n\to \infty} \frac{\#\{0\le i \le n: \widetilde{T}^i(x,M)\in (C_s,M')\}}{n} = \mu (C_s,M').
\]

We say $\widetilde{T}$ is transitive if for any $M_1,M_2\in \mathcal{M}$, there exists a string $s$ such that \[ T^n( C_s,M_1) = (\Omega,M_2), \] where $n$ is the length of the string $s$. This string is called the traversing string from $M_1$ to $M_2$.

\begin{thm}\label{thm:traversing}
If $\widetilde{T}$ is transitive and is measure-preserving with respect to $\mu$, then  $\widetilde{T}$ is ergodic.  Moreover, if $x$ is normal, then for any $M\in \mathcal{M}$, the point $(x, M)$ is $\widetilde{T}$-normal with respect to $\mu$.
\end{thm}

In \cite{ratmultCF}, this result was proved for the continued fraction Gauss map. However, the only facts we used about the Gauss map are also satisfied for all the base-$b$ maps, see Theorem 3.1 and Remark 3.2 in \cite{ratmultCF}. Thus this theorem also holds in this case. 

The usefulness of Theorem \ref{thm:traversing} is that it tells us that any normal number lifted to the augmented system will still be normal in this augmented system.

\section{Proof of Theorem \ref{thm:leap}}

Let $k$ be a fixed positive integer. Consider the augmented system given by $\mathcal{M}=(\ell, [b_1,b_2,\dots,b_k])$ with $\ell, b_1,b_2,\dots,b_k \in \{0,1,2,\dots,b-1\}$, $\ell\le b_k$, and 
\[
\widetilde{T}(x,\ell,[b_1,b_2,\dots,b_k]) = \begin{cases}
(Tx, \ell-1, [b_1,b_2,\dots,b_k]), & \text{if }\ell\ge 1,\\
(Tx, a_1(x), [b_2,b_3,\dots,b_k,a_1(x)]), & \text{if } \ell=0.
\end{cases}
\] (We suppress extra parentheses for readability.) 
For our measure, we just use $\mu=\lambda \times c$ where $\lambda$ is the usual Lebesgue measure on $[0,1)$ and $c$ is the normalized counting measure on $\mathcal{M}$. 

To apply Theorem \ref{thm:traversing}, we must show that $\widetilde{T}$ is both traversing and measure-preserving. 

For transitivity, let $M_1=(\ell,[b_1,b_2,\dots,b_k])$ and $M_2 = (\ell', [b_1',b_2',b_3',\dots,b_k'])$. Let $1^j$ denote $j$ repetitions of the digit $1$. Then the desired traversing string from $M_1$ to $M_2$ is given by
\[
[1^\ell,b_1',1^{b_1'},b_2',1^{b_2'},b_3',\dots,1^{b_{k-1}'},b_k',1^{b_k'-\ell'}].
\]

For measure-preserving, consider the inverse branches $T_j^{-1}$ of $T$ that are defined by $T_j^{-1} 0.a_1a_2\dots = 0.ja_1a_2\dots$. We define the inverse branches $\widetilde{T}_j^{-1}$ of $\widetilde{T}$ by being the branch of the inverse map that induces $T_j^{-1}$ in the first coordinate. For any measurable set $E\subset \Omega$,
\begin{equation}\label{eq:firstinverse}
\widetilde{T}_j^{-1}(E,\ell,[b_1,b_2,\dots,b_k]) = \begin{cases}
(T_j^{-1} E, \ell+1, [b_1,b_2,\dots,b_k]), & \text{if }\ell< b_k,\\
(T_j^{-1} E, 0, [*, b_1,b_2,\dots,b_{k-1}]), & \text{if } j=\ell=b_k,\\
\emptyset , & \text{otherwise.}\end{cases}
\end{equation}
In this paper we use $*$ to refer to an arbitrary digit, so the second case of \eqref{eq:firstinverse} is actually a union over sets of the form $(T_j^{-1} E, 0, [d, b_1,b_2,\dots,b_{k-1}])$ as $d$ runs over all elements of the set $\{0,1,\dots,b-1\}$.

The measure of the set $(E,\ell, [b_1,b_2,\dots,b_k])$ is $\lambda(E)/|\mathcal{M}|$. On the right-hand side of \eqref{eq:firstinverse}, the measure is $\lambda(E)/b\cdot|\mathcal{M}|$ in the first case, $ \lambda(E)/|\mathcal{M}|$ in the second case, and $0$ in the last case. By summing over $j$, we see that $\widetilde{T}$ preserves the measure of sets of the form $(E,\ell,[b_1,b_2,\dots,b_k])$. 

Now consider a more general set $\widetilde{E}\subset\widetilde{\Omega}$. We may write $\widetilde{E}$ as a disjoint union of sets of the form $(E_M,M)$ for $M\in\mathcal{M}$. If the inverse images $\widetilde{T}^{-1}(E_M,M)$ are all disjoint, then it is clear that $\widetilde{T}$ preserves $\mu$. Consider $M=(\ell, [b_1,b_2,\dots,b_k])$ and $M'=(\ell',[b'_1,b'_2,\dots,b'_k])$ with $M\neq M'$. If $\ell< b_k$ and $\ell'< b'_k$, then the action of $\widetilde{T}^{-1}$ on the $\mathcal{M}$-coordinate is just to increase $\ell$ and $\ell'$ by one, and thus  $\widetilde{T}^{-1}(E_M,M)$ and $\widetilde{T}^{-1}(E_{M'},M')$ must still differ in the $\mathcal{M}$-coordinate. This is still true if $\ell < b_k$ but $\ell' = b'_k$, since in this case, $\ell$ is taken to $\ell+1$ but $\ell'$ is taken to $0$, and these are clearly different values. It remains to consider the case where both $\ell = b_k$ and $\ell'=b'_k$. In this case, we have
\[
\widetilde{T}^{-1}(E_M,M) = (T_{b_k}^{-1}E_M, 0,[*,b_1,b_2,\dots,b_{k-1}]),
\]
so if $b_i\neq b'_i$ for some $i$ with $1\le i \le k-1$, we clearly have that $\widetilde{T}^{-1}(E_M,M)$ and $\widetilde{T}^{-1}(E_{M'},M')$ differ in the $\mathcal{M}$-coordinate; otherwise, we must have that $b_k\neq b'_k$, and then we get the same result in the first coordinate since the inverse branches $T^{-1}_{b_k}$ and $T^{-1}_{b'_k}$ have disjoint images.

Thus, we have that $\widetilde{T}$ preserves $\mu$, and Theorem \ref{thm:traversing} applies.

Now we return to the points $x=0.a_1a_2a_3\dots$ that we assume is normal to base $b$ and the point $y=0.a_{n_1}a_{n_2}a_{n_3}\dots$ given by $n_{i+1} = n_i+1+a_{n_i}$, $ i\ge 1$, which we want to prove is normal to base $b$. We may assume without loss of generality that $n_1=1$; if it were larger, we could replace $x$ with $0.a_{n_1}a_{n_1+1}a_{n_1+2}\dots$, which is also normal to base $b$. Consider a string $s=[c_1,c_2,\dots,c_k]$ and consider the augmented system for this $k$ as detailed in the previous paragraphs. Let $\tilde{x}=(x,0,[0,0,\dots,0])$ in the augmented space, and let $\mathcal{N}$ denote the subset of $\mathcal{M}$ such that the corresponding $\ell$ and $b_k$ are equal.  By Theorem \ref{thm:traversing}, we know that $\tilde{x}$ is $\widetilde{T}$-normal with respect to $\mu$. 

This augmented system was specifically constructed so that \[ \widetilde{T}^{n_i} \tilde{x} = (T^{n_i} x, a_{n_i},[ a_{n_{i-k+1}}, a_{n_{i-k+2}}, \dots, a_{n_i}])\] where we are assuming that $a_{j}=0$ if $j\le 0$. In fact, we have $\widetilde{T}^n\tilde{x} \in [0,1)\times \mathcal{N}$ with $n>0$ if and only if $n=n_i$ for some $i$. Therefore if $T^i y \in C_s$, then $\widetilde{T}^{n_{i+k}} \tilde{x} \in ([0,1),c_k,s)$, and these are the only times the $\widetilde{T}^n \tilde{x}$ visits this set (provided $n>n_k$).

 Thus,
\[
 \frac{\#\{0\le i \le m-1: T^i y \in C_s\}}{m}=  \frac{ \#\{0\le i \le n: \widetilde{T}^i \tilde{x} \in  ([0,1), c_k,s)\}+O(1)}{ \# \{0\le i \le n:  \widetilde{T}^i \tilde{x} \in[0,1)\times \mathcal{N}\}+O(1)}
\]
for any $n$ with $n_{m+k}\le n < n_{m+k+1}$, provided $m$ is sufficiently large. Here, the $O(1)$ terms are to account for irregularities caused for small values of $i$. Therefore, we have
\begin{align*}
\lim_{m\to \infty} \frac{\#\{0\le i \le m-1: T^i y \in C_s\}}{m}  &= \lim_{n\to \infty} \frac{ \#\{0\le i \le n: \widetilde{T}^i \tilde{x} \in  ([0,1), c_k,s)\}}{ \# \{0\le i \le n:  \widetilde{T}^i \tilde{x} \in[0,1)\times \mathcal{N}\}}\\ &= \lim_{n\to \infty} \frac{ \#\{0\le i \le n:  \widetilde{T}^i \tilde{x} \in  ([0,1), c_k,s)\} /n }{ \# \{0\le i \le n:  \widetilde{T}^i \tilde{x} \in [0,1)\times \mathcal{N}\} /n }\\
&= \frac{\mu([0,1),c_k,s)}{\mu([0,1)\times \mathcal{N})}=\frac{1/|\mathcal{M}|}{|\mathcal{N}|/|\mathcal{M}|}=\frac{1}{|\mathcal{N}|},
\end{align*}
where the  equality between the last two lines follows by applying the definition of $\tilde{x}$ being $\widetilde{T}$-normal with respect to $\mu$.  Finally, it is easy to count that $|\mathcal{N}|=b^{-k}$, so that the string $s$ does occur in the base-$b$ expansion of $y$ with the desired frequency.

\section{Proof of Theorem \ref{thm:reduce}}

The proof here (and in subsequent proofs) is so similar to the proof of Theorem \ref{thm:leap} that we will restrict ourselves to proving that a nice augmented system with the desired properties exists, rather than repeating identical steps of the remainder of the proof.

In this case we consider an augmented system (for an arbitrary fixed, positive integer $k$) that is given by $\mathcal{M}=([b_1,b_2,\dots,b_k])$, $b_i \in \{0,1,2,\dots,b-2\}$,  $\widetilde{T}$ given by
\[\widetilde{T}(x,[b_1,b_2,\dots,b_k]) = \begin{cases}
(Tx, [b_2,b_3,\dots,b_k, a_1(x)]), & \text{if } a_1(x) < b-1\\
(Tx, [b_1,b_2,\dots,b_k ]), & \text{otherwise},
\end{cases}
\]
and we take $\mu$ again to be Lebesgue measure crossed with the normalized counting measure on $\mathcal{M}$.

Here, transitivity is  trivial. If $M_1 = ([b_1,b_2,\dots,b_k])$ and $M_2 = ([b'_1,b'_2,\dots,b'_k])$, then the traversing string from $M_1$ to $M_2$ is given by $[b'_1,b'_2,\dots,b'_k]$ itself. 

Now we want to show measure-preserving. Consider again the inverse branches $\widetilde{T}_j^{-1}$ of the transformation $\widetilde{T}$. For any measurable set $E\subset\Omega$, we have
\begin{equation}\label{eq:secondinverse}
\widetilde{T}_j^{-1} (E,[b_1,b_2,\dots,b_k]) = \begin{cases}
(T_j^{-1} E, [b_1,b_2,\dots,b_k]), & \text{if } j=b-1\\
(T_j^{-1} E, [*, b_1, b_2,\dots,b_{k-1}]), & \text{if } j=b_k\\
\emptyset, & \text{otherwise}.
\end{cases}
\end{equation}
The measure of the set $(E,[b_1,b_2,\dots,b_k])$ is $\lambda(E)/(b-1)^k$. The measure of the set  on the right-hand side of \eqref{eq:secondinverse} is $\lambda(E)/b(b-1)^k$ in the first case, $\lambda(E)/b(b-1)^{k-1}$ in the second case, and $0$ in the third case. Summing over $j$ we see that $\widetilde{T}$ preserves the measure of sets of this form.

So we consider an arbitrary set $\widetilde{E}\subset\widetilde{\Omega}$ and write it as the disjoint union of sets of the form $(E_M,M)$. Let $(E_M,M)$ and $(E_{M'},M')$ be two such sets and consider the inverse image of both of them. We will show these must be distinct if $M$ and $M'$ are distinct. Clearly $T_j^{-1} (E_M,M)$ and $T_{j'}^{-1}(E_{M'},M')$ are disjoint if $j\neq j'$, as they will be disjoint in the first coordinate. So assume that $j=j'$ for a moment. If $j=b-1$, then the only way for  $T_j^{-1} (E_M,M)$ and $T_{j'}^{-1}(E_{M'},M')$ to overlap is if $M=M'$, which would go against our assumption that they are distinct. If $j<b-1$ and $j \neq b_k$ (or $j \neq b'_k$), then $T_j^{-1} (E_M,M)=\emptyset$  (or $T_j^{-1}(E_{M'},M')= \emptyset$) and the result is trivial in this case. Finally if $j<b-1$ and $j=b_k=b'_k$, then it is easy to see that the only way for  $T_j^{-1} (E_M,M)$ and $T_{j'}^{-1}(E_{M'},M')$ is if we have $b_i=b'_i$ for $1\le i \le k-1$, so we can again see that they are distinct unless $M=M'$, which proves that $\widetilde{T}$ is measure-preserving. Thus Theorem \ref{thm:traversing} applies.

We finish by remarking on how to count the necessary occurrences of strings and $n_i$'s. In this case, we let $\tilde{x}= (x,[0,0,\dots,0])$ and let $y=0.a_{n_1}a_{n_2}a_{n_3}\cdots$ as before. Then $T^i y \in C_s$ if and only if $\widetilde{T}^{n_{i+k+1}-1} \tilde{x} \in ([0,1-1/b), [c_1,c_2,\dots,c_k])$ for sufficiently large $i$. Likewise we have that $n=n_i$ if and only if $\widetilde{T}^{n-1} \tilde{x} \in ([0,1-1/b), *)$ for sufficiently large $n$. The  set $[0,1-1/b)$ is important because it is the set of $x$'s with $a_1(x)\neq b-1$, and thus represents the set on which if $\tilde{T}$ is applied, then the digits of $\mathcal{M}$ are shifted one position.

Since we have
\[
\frac{ \mu ([0,1-1/b), [c_1,c_2,\dots,c_k])}{\mu ([0,1-1/b), *)} = \frac{ \frac{b-1}{b}\cdot \frac{1}{(b-1)^k}}{\frac{ b-1}{b}} = \frac{1}{(b-1)^k} = \lambda(C_s),
\]
where here the cylinder set $C_s$ is thought of as a cylinder set of length $k$ for the base-$(b-1)$ transformation, the theorem holds.

\section{Proof of Theorem \ref{thm:modulo}}

In this case, consider an augmented system (for an arbitrary fixed, positive integer $k$) that is given by $\mathcal{M}=(\ell,[b_1,b_2,\dots,b_k])$ with $\ell \in \{0,1,2,\dots,N-1\}$, $b_i\in \{0,1,\dots, b-1\}$, $\widetilde{T}$ is given by
\begin{align*}
&\widetilde{T}(x,\ell,[b_1,b_2,\dots,b_k]) \\ &\quad= \begin{cases}
(Tx, \ell+a_1(x) \ (\bmod{N}), [b_1,b_2,\dots,b_k]), & \text{if }\ell+a_1(x)\not\equiv L(\bmod{N}),\\
(Tx, \ell+a_1(x)\  (\bmod{N}), [b_2,b_3,\dots,b_k,a_1(x)], & \text{if }\ell +a_1(x)\equiv L (\bmod{N}),
\end{cases},
\end{align*}
where, in the second coordinate, we are assuming all values are reduced to be in the set $\{0,1,2,\dots, N-1\}$, and finally we take $\mu$ again to be Lebesgue measure crossed with the normalized counting measure on $\mathcal{M}$. We must show $\widetilde{T}$ is transitive and measure-preserving.

For transitivity, again  let $M_1=(\ell,[b_1,b_2,\dots,b_k])$ and $M_2 = (\ell', [b_1',b_2',b_3',\dots,b_k'])$. Then the desired string is given by
\begin{align*}
&[1^{\ell-b_1'\bmod{N}},b_1',1^{-b_2'\bmod{N}},b_2',1^{-b_3'\bmod{N}},b_3',\\ &\qquad\dots,1^{-b_{k}'\bmod{N}},b_k',1^{\ell'-L\bmod{N}}].
\end{align*}

Finally it remains to show measure-preserving. Again we consider inverse branches $\widetilde{T}_j^{-1}$ of the transformation $\widetilde{T}$. For any measurable set $E\subset \Omega$, we have
\begin{align*}
&\widetilde{T}_j^{-1}(E,\ell,[b_1,b_2,\dots,b_k]) \\ &\quad= \begin{cases}
(T_j^{-1} E, \ell-j\ (\bmod{N}), [b_1,b_2,\dots,b_k]), & \text{if }\ell\neq L,\\
(T_j^{-1} E, \ell-j\ (\bmod{N}), [*, b_1,b_2,\dots,b_{k-1}]), & \text{if } \ell=L, \quad j=b_k,\\
\emptyset , & \text{otherwise.}\end{cases}
\end{align*}
The measure of $(E,\ell, [b_1,b_2,\dots,b_k])$ is $\lambda(E)/|\mathcal{M}|$. The measure of the set on the right-hand side above is also $\lambda(E)/b\cdot |\mathcal{M}|$ in the first case, $\lambda(E)/|\mathcal{M}|$ in the second case, and $0$ in the third case. By summing over $j$, we see that $\widetilde{T}$ preserves the measure of sets of the form $(E,M)$. 

Consider a more general set $\widetilde{E}\subset\widetilde{\Omega}$. We may write $\widetilde{E}$ as a disjoint union of sets of the form $(E_M,M)$ for $M\in\mathcal{M}$. Again we want to show that the sets $\widetilde{T}^{-1}_j (E_M,M)$ and $\widetilde{T}^{-1}_{j'} (E_{M'},M')$ are disjoint unless $j=j'$ and $M=M'$. Let $M=(\ell, [b_1,b_2,\dots,b_k])$ and $M' = (\ell' , [b'_1,b'_2,\dots,b'_k])$. Clearly these sets are disjoint in the first coordinate if $j\neq j'$, so suppose $j=j'$. If $\ell \neq \ell'$ then since $j=j'$, the two sets would be disjoint in the second coordinate, so we suppose that $j=j'$ and $\ell=\ell'$. But then it is easy to see that if $\ell = \ell'$ and $j=j'$ but $M \neq M'$, then it is clear that the two sets will be disjoint, which completes the final case needed to show that $\widetilde{T}$ is measure-preserving. Thus Theorem \ref{thm:traversing} applies.

We finish again by remarking on how to count the necessary occurrences of strings and $n_i$'s. In this case, we let $\tilde{x}= (x,0,[0,0,\dots,0])$ and let $y=0.a_{n_1}a_{n_2}a_{n_3}\cdots$, and then $T^i y \in C_s$ if and only if $\widetilde{T}^{n_{i+k}} \tilde{x} \in (\Omega ,L, [c_1,c_2,\dots,c_k])$, for $i$ sufficiently large. Likewise we have that $n=n_i$ if and only if $\widetilde{T}^{n} \tilde{x} \in (\Omega, L, *)$, for $n$ sufficiently large.

Since we have
\[
\frac{ \mu (\Omega ,L, [c_1,c_2,\dots,c_k])}{\mu(\Omega,L,*)} = \frac{\frac{1}{N}\cdot \frac{1}{b^k}}{\frac{1}{N}} = \frac{1}{b^k} = \mu(C_s),
\]
the theorem holds.

\section{Acknowledgments}

The author acknowledges assistance from the Research and Training Group grant DMS-1344994 funded by the National Science Foundation.

The author would also like to thank Andy Parrish for his help and insights.

\end{document}